\documentclass[12pt]{amsart}
\usepackage{amssymb,amsmath,amsthm,amsfonts,times,hyperref,mathrsfs,multirow}
\newtheorem{theorem}{Theorem}[section]
\newtheorem{lemma}[theorem]{Lemma}
\newtheorem{cor}[theorem]{Corollary}

\theoremstyle{definition}

\newtheorem{example}[theorem]{Example}

\theoremstyle{remark}
\newtheorem*{remark}{Remark}

\numberwithin{equation}{section}

\newcommand{\SLZ}{\mbox{SL}_2(\mathbb{Z})}

\newcommand{\Parans}[1]{\left(#1\right)}

\newcommand\leg[2]{\genfrac{(}{)}{}{}{#1}{#2}} 

\newcommand\twidit[1]{\overset {\text{\lower 3pt\hbox{$\sim$}}}#1}
\newcommand\dtwidit[1]{\overset {\text{\lower 6pt\hbox{$\sim$}}}#1}
\newcommand\Wtwid{\overset {\text{\lower 3pt\hbox{$\sim$}}}W}
\newcommand\gtwid{\overset {\text{\lower 3pt\hbox{$\sim$}}}g}
\newcommand\ttwid{\overset {\text{\lower 3pt\hbox{$\sim$}}}\theta}
\newcommand\mutwid{\overset {\text{\lower 3pt\hbox{$\sim$}}}\mu}

\newcommand\AMat{\begin{pmatrix} a & b \\ c & d \end{pmatrix}}

\newcommand\EObara{\overline{\mathcal{EO}}}
\newcommand\EObar[1]{\overline{\mathcal{EO}}({#1})}
\allowdisplaybreaks
\newcommand\mylabel[1]{\label{#1}}
\newcommand\thm[1]{\ref{thm:#1}}
\newcommand\omylabel[1]{\quad\framebox{#1}\label{#1}}
\newcommand\myeqn[1]{(\ref{eq:#1})\framebox{#1}}

\newcommand\omycite[1]{}

\newcommand\eqn[1]{(\ref{eq:#1})}
\newcommand\sect[1]{\ref{sec:#1}}

\newcommand\subsect[1]{\ref{subsec:#1}}

\newcommand{\beqs}{\begin{equation*}}
\newcommand{\eeqs}{\end{equation*}}
\newcommand{\beq}{\begin{equation}}
\newcommand{\eeq}{\end{equation}}

\newcommand{\sB}{\mathscr{B}}  

\newcommand\Tr[3]{\mbox{\rm Tr}({#1}, {#2})( {#3} )}

\newcommand{\PF}[2]{\Parans{\frac{#1}{#2}}}
\newcommand{\pp}{\overline{p}}
\newcommand{\qq}{\overline{q}}
\newcommand{\podd}{\mbox{pod}}

\begin{document}
\renewcommand{\MR}[1]{\href{http://www.ams.org/mathscinet-getitem?mr={#1}}{MR{#1}}}
\title[Multiplicative Congruences]
{Multiplicative Congruences for Andrews's Even Parts Below Odd
Parts Function and Related Infinite Products}
\author{F. G. Garvan}
\address{Department of Mathematics, University of Florida,
Gainesville,
FL 32611-8105}
\email{fgarvan@ufl.edu}
\author{Connor Morrow}
\address{Department of Mathematics, University of Florida,
Gainesville,
FL 32611-8105}
\email{connorfmorrow@ufl.edu}

\subjclass[2020]{11P82,11P83,11F37}

\date{\today}

\dedicatory{Dedicated to our friend and colleague, Mourad Ismail, on the occasion of his 80th birthday}

\keywords{Partition congruences, eta-quotients, Hecke operators, crank parity function, Andrews $\EObara$-function}

\begin{abstract}
We prove multiplicative congruences mod $2^{12}$ for George
Andrews's partition function, $\overline{\mathcal{EO}}(n)$, the
number of partitions of $n$  in which every even part is less than
each odd part and only the largest even part occurs an odd number
of times. We find analogous congruences for more general infinite
products. These congruences are obtained using
Fricke involutions and Newman's approach to half integer weight Hecke
operators on eta quotients, and were inspired by Atkin's multiplicative
congruences for the partition function.
\end{abstract}

\maketitle
\section{Introduction}

Andrews \cite{An2018} studied the
partition function $\mathcal{EO}(n)$ which counts the number of
partitions of $n$ where every even part is less than each odd
part. He denoted by $\overline{\mathcal{EO}}(n)$, the number of
partitions counted by $\mathcal{EO}(n)$ in which \textit{only}
the largest even part appears an odd number of times. For example,
$\mathcal{EO}(8)=12$ with the relevant partitions being
\begin{gather*}
8,
6+2, 7+1, 4+4, 4+2+2, 5+3, 5+1+1+1, 2+2+2+2, 3+3+2,\\
 3+3+1+1, 3+1+1+1+1+1, 1+1+1+1+1+1+1+1+1;
\end{gather*}
and $\overline{\mathcal{EO}}(8)=5$,
with the relevant partitions being
$$
8, 4+2+2, 3+3+2, 3+3+1+1,
1+1+1+1+1+1+1+1.
$$

We use the standard $q$-series notation:
     \begin{equation*}
         (a; q)_n = \prod_{k=0}^{n-1}(1-aq^k),
     \end{equation*}
\begin{equation*}
    (a;q)_\infty = \lim_{n \rightarrow \infty} (a; q)_n =
    \prod_{k=0}^{\infty}(1-aq^k) \hspace{5mm} (|a| < 1),
\end{equation*}
and
\begin{equation}
    E(q) = (q)_{\infty} = (q;q)_{\infty}.
    \omylabel{eq:eqdef}
\end{equation}

Andrews found
   \begin{equation*}
\sum_{n=0}^{\infty}\overline{\mathcal{EO}}(n)q^n
=\prod_{n=1}^\infty
\frac{ (1-q^{4n})^3 }{ (1-q^{2n})^2},
\end{equation*}
and proved the following congruence:
   \begin{equation*}
    \overline{\mathcal{EO}}({10n+8}) \equiv 0 \, (\bmod{\,5}).
\end{equation*}

He also found a crank-type function to explain this
congruence. Previously, $\overline{\mathcal{EO}}(n)$
 congruences had been found mod 4 by Chen and Chen \cite{Ch-Ch2023},
 Ray and Barman mod 8 \cite{Ra-Ba2020}, and Chen \cite{ChSC2023}
 mod 2, 4, and 8. We extend the congruences to mod $2^{12}$, and some
explicit congruences mod $2^i$ for $1 \le i \le 5$.
 To simplify results, we let
\begin{equation}
\mylabel{eq:alhadef}
  \sum_{n=0}^\infty \alpha(n) q^n =
\prod_{n=1}^\infty
\frac{ (1-q^{2n})^3 }{ (1-q^{n})^2},
\end{equation}
so that
\begin{equation*}
    \overline{\mathcal{EO}}(n) =
\begin{cases}
\alpha(m) & \mbox{if $n=2m$}\\
0 & \mbox{otherwise.}
\end{cases}
\end{equation*}
We mainly consider congruences for $\alpha(n)$ but also
consider analogous infinite products.
One of our main results is
\begin{theorem}
\mylabel{thm:alcong212}
Let $\ell > 3$ be prime, and $N_{\ell} = \frac{1}{6}(\ell^2-1).$ Then
\begin{equation}
\mylabel{eq:212cong}
     \ell \alpha(\ell^2n + N_{\ell}) +
     \left(\frac{-1}{\ell}\right)\left(\frac{N_{\ell}-n}{\ell}\right)\alpha(n)
         + \alpha\left(\frac{n-N_{\ell}}{\ell^2}\right) \equiv
         c_{\ell} \, \alpha(n) \, \, (\bmod{\,2^{12}})
\end{equation}
    where $c_{\ell} =
    \ell\alpha(N_{\ell})+\left(\frac{6}{\ell}\right)$,
    and $\left(\frac{\cdot}{\ell}\right)$ is the Legendre symbol.
\end{theorem}

\begin{example}
For $\ell = 5$ we find that $c_5=26$ so that
\begin{equation*}
    5\alpha(25n+4) + \left(\frac{4-n}{5}\right) \alpha(n) +
    \alpha\left(\frac{n-4}{25}\right) \equiv 26 \alpha(n) \, \,
    (\bmod{\,2^{12}}).
\end{equation*}
\end{example}
\begin{example} For $\ell = 31$ we find that
    \begin{equation*}
        c_{31} = 2^{10}\cdot 5^2 \cdot 222823,
    \end{equation*}
    and it follows that
    \begin{equation*}
        \alpha(31^3n+4965)\equiv 0 \, \, (\bmod{\,2^{10}})
    \end{equation*}
    provided $n \not \equiv 5 \, \, (\bmod{\,31})$.
\end{example}
By a careful study of the constant, $c_{\ell},$ we prove
\begin{theorem}
\mylabel{thm:alcongs32}
Let $1 \le i \le 5$ and assume $3 \ne \ell \equiv -1 \pmod{2^i}$ is
prime. Then
\begin{equation}
\mylabel{eq:232cong}
     \alpha(\ell^2(\ell n+m)+N_{\ell}) \equiv 0 \, \,
     (\bmod{\,2^{i}})
\end{equation}
where $n \not\equiv N_{\ell} \, \, (\bmod{\,\ell})$, and $m>0$
is such that $m \equiv N_{\ell} \, \, (\bmod{\,\ell})$.
\end{theorem}

Our proof of Theorem \thm{alcong212} depends on some results of
Newman \cite{Ne1962} for half-integer weight Hecke operators acting on 
eta-quotients. 
We note that Theorems \thm{alcong212} and \thm{alcongs32}
were discovered independently  by the second author in his thesis
\cite{Mo-thesis}, where full details are given.
Newman's results enable us to consider more
congruences for more general eta-quotients.  We consider
eta-quotients of the form
\beq
\mylabel{eq:crspdef}
\sum_{n=0}^\infty c_{r,s,p}(n) q^n = E(q)^r \, E(q^p)^s,
\eeq
where
$$
E(q) = (q)_\infty = q^{-1/24} \eta(\tau).
$$
We obtain the following generalization of Theorem \thm{alcong212}. 
\begin{theorem}
\mylabel{thm:gencong}
Suppose $p=2$, $3$, $5$, $7$ or $13$ and $3<\ell\ne p$ is prime.
Suppose that $s>0$ and $r<0$ are integers of opposite parity satisfying
$$
0\le r + sp < 24 \qquad \mbox{and} \qquad s + rp < 0, 
$$ and define
$$
c(n) = c_{r,s,p}(n).
$$
Then
\begin{equation}
\mylabel{eq:gencong}
\ell^{2-2\varepsilon} c(n\ell^2+\Delta)
+ \left(\frac{\theta}{\ell}\right)
\ell^{\frac{1}{2}-\varepsilon}\left(\frac{n-\Delta}{\ell}\right)c(n)+
c\left(\frac{n-\Delta}{\ell^2}\right) \equiv \lambda_{r,s,p,\ell}\,
c(n) \, \, (\bmod{\,p^{\frac{k}{2}}})
\end{equation}
where 
$$
\theta = (-1)^{\frac{1}{2}-\frac{1}{2}({r}+s)}2p^{s}, \quad
\Delta = (r+sp)\frac{(\ell^2-1)}{24},\quad  k = \frac{24}{p-1},\quad
\varepsilon = \frac{1}{2}(r + s),
$$
and
\begin{equation*}
    \lambda_{r,s,p,\ell} = \ell^{2-2\varepsilon}
    c(\Delta) + \left(\frac{\theta}{\ell}\right)
    \ell^{\frac{1}{2}-\varepsilon}\left(\frac{-\Delta}{\ell}\right)+
    c\left(\frac{-\Delta}{\ell^2}\right).
\end{equation*}
\end{theorem}

In Section \sect{modfuncs+newman} we review the definition of
modular function and summarize Newman's results on eta-quotients
and his approach to half-integer weight Hecke operators.
In Section \sect{proofs}  we apply Newman's results to prove 
our main result, Theorem \thm{gencong} which implies Theorem \thm{alcong212},
our mod $2^{12}$ congruence for the $\EObara$-function.
In Section \sect{fabpoly} we review Faber polynomials and their connection
to the $j$-function, and recall Ono's analog and how it
relates to Atkin's multiplicative congruences for the partition function.
We obtain analogues for more general eta-quotients corresponding to 
genus zero cases, and further congruences for the eigenvalues in Theorem
\thm{gencong}, which also leads to Theorem \thm{alcongs32}.
Finally in Section \sect{end} we relate the crank parity and 
$\EObara$-function to identities for Ramanujan's third order
mock theta functions, and we pose some problems.

\section{Modular Functions and Newman's Results}
\mylabel{sec:modfuncs+newman}

\subsection{Modular Functions}
\mylabel{subsec:modfuncs}
For $\gamma = \AMat\in \SLZ$ we define
$$
\gamma \tau := \frac{a \tau + b}{c \tau +d}
$$
for $\tau$ in the complex upper-half plane $\mathscr{H}$. 
A subgroup $\Gamma'$ of $\SLZ$ has \textit{level}  $N$ if it contains
\begin{equation*}
 \Gamma(N) = \{\gamma \in \SLZ\,:\,  \gamma \equiv I \,
        (\bmod{\,N}) \}.
\end{equation*}
We say a function $\,:\,\mathscr{H} \longrightarrow \mathbb{C}$ is 
a \textit{modular function} on a subgroup $\Gamma'$ of $\SLZ$ of level
$N$ if its satisfies

\begin{enumerate}
\item[(i)]
$f$ is meromorphic on $\mathscr{H}$;
\item[(ii)]
$$
f(A \tau) = f(\tau),
$$
for $A\in \Gamma'$; and
\item[(iii)]
For all $\gamma \in \SLZ$, then $f(\gamma \tau)$ 
 has a Fourier expansion of the form
\begin{equation*}
        f(\gamma \tau) = \sum_{n\geq n_0} a_{\gamma}(n)q_N^n,
    \end{equation*}
    where $q_N := e^{2\pi i \tau / N}$ and $a_{\gamma}(n_0) \neq 0.$
\end{enumerate}
If in (iii), $n_0<0$, then we say $f$ has a pole at the cusp
$\gamma(i\infty)=a/c$. For brevity we say $f$ \textit{is on} $\Gamma'$. 
We also say $f$ is \textit{weakly holomorphic} if the only
poles are at cusps.

We will mainly be consider modular functions on
$$
\Gamma_0(N) = \left\{ \AMat \in \SLZ\,:\, c \equiv 0 \pmod{N}\right\}.
$$
\subsection{Newman's Results}
\mylabel{subsec:newmanres}
Throughout we assume $r$, $s$ are integers and $p$ and $\ell$
are distinct primes.
If $f = f(r,s)$, then we shall write $f^* = f(s,r)$.
We define
\beq
\mylabel{eq:Bdef}
B(\tau) = B_{r,s,p}(\tau) = \eta^r(\tau) \, \eta^s(p \tau),
\eeq
where
$$
\eta(\tau) = q^{1/24} \prod_{n=1}^\infty (1 - q^n),
$$
and we define
\beq
\mylabel{eq:phidef}
\phi(\tau) = \phi_{r,s,p}(\tau) = 
\prod_{n=1}^\infty (1 - q^n)^r 
\prod_{n=1}^\infty (1 - q^{pn})^s 
=
\sum_{n=0}^\infty c_{r,s,p}(n) q^n = 
\sum_{n=0}^\infty c(n) q^n.       
\eeq
We let
$$
\varepsilon = \frac{1}{2}(r + s), \quad t = \frac{1}{24}(r + sp), 
$$
so that 
\begin{equation*}
    B(\tau) = q^t \phi(\tau).
\end{equation*}
We have
\begin{equation*}
    B^*(\tau) = \eta^{r}(p\tau) \eta^s(\tau) = q^{t^*} \phi^*(\tau),
\end{equation*}
where
\beq
\mylabel{eq:phistardef}
    \phi^*(\tau) = \prod_{n=0}^{\infty}(1-q^{pn})^{r}
    \prod_{n=0}^{\infty}(1-q^{n})^s 
    =\sum_{n=0}^\infty c^*(n) q^n. \qquad
    t^* = \frac{1}{24}({r}p+s).
\eeq
We define 
$$
g(\tau) = g_{r,s,p,\ell}(\tau) = \frac{B(\ell^2 \tau)}{B(\tau)}.
$$
We have
\begin{lemma}[{ Newman \cite[Lemma 2]{Ne1962} }]
\mylabel{lem:Nelem2}
If $\Delta= t(\ell^2-1) = ({r}+sp)\frac{(\ell^2-1)}{24}$ and
$\Delta^{*} =  t^*(\ell^2-1) = (s+rp)\frac{(\ell^2-1)}{24}$
are integers then $g(\tau)$ is a weakly holomorphic modular
functions on
$\Gamma_0(\ell^2 p)$.
\end{lemma}
We remark that if $\ell > 3$ is prime then $\frac{1}{24}(\ell^2-1)$,
$\Delta$, and $\Delta^{*}$ are integers.

We describe Newman's method of subgroups. 
Instead of  half-integer weight Hecke operators, Newman uses trace maps.
Let $\Gamma_0 < \Gamma_1$ be subgroups of
finite index in $\SLZ$. Suppose the function $f(\tau)$ is on $\Gamma_0$
then the function
$$
\Tr{\Gamma_0}{\Gamma_1}{f} := \sum_{j=1}^n f(R_n \tau)
$$
is a function on the group $\Gamma_1$, 
where $R_1$, $R_2$, \dots, $R_n$ are right coset representatives of $\Gamma_0$
in $\Gamma_1$. See Newman \cite[Theorem 2.2]{Ne1952}.     

Let $W =  \begin{pmatrix}
1 & 0\\
1 & 1
\end{pmatrix}$. Then, $W^{-k\ell p}$ with $0\leq k \leq \ell-1$
forms a set of right coset representatives for $\Gamma_0(\ell^2 p)$
in $\Gamma_0(\ell p)$, and
\begin{equation*}
    h(\tau) = \sum_{k=0}^{\ell-1}g(W^{-k\ell p}\tau) =
    \Tr{\Gamma_0(\ell^2 p)}{\Gamma_0(\ell p)}{g}
\end{equation*}
is a modular function on $\Gamma_0(\ell p)$. We let
$R = \begin{pmatrix}
        1 && -1 \\
        -pp_0 && \ell_0\ell
    \end{pmatrix}$  where $\ell_0\ell - p_0p = 1$, and $W^{-kp}$ with
    $0\leq k \leq \ell-1$ form a set of right coset representatives
    for $\Gamma_0(\ell p)$ in $\Gamma_0(p),$ and
    \begin{equation*}
\Tr{\Gamma_0(\ell p)}{\Gamma_0(p)}{h} = 
        \sum_{k=0}^{\ell}h(W^{-kp}\tau) + h(R\tau)
    \end{equation*}
is a modular function on $\Gamma_0(p)$. This leads to the following
theorem.
\begin{theorem}[{ Newman \cite[Lemma 3]{Ne1962} }] 
\mylabel{thm:Gthm}
Define \begin{equation*}
        G(\tau) = G_{r,s,\ell, p} (\tau) = 
         \sum_{k=0}^{\ell^2-1}g(W^{-kp}\tau) +
        \sum_{k=0}^{\ell-1}g(W^{-k\ell p}R\tau).
    \end{equation*}
    Then, $G(\tau)$ is a weakly holomorphic modular function on
    $\Gamma_0(p)$. Also for $2\varepsilon$ odd,
\begin{equation*}
    \phi(\tau) \, G(\tau) = 
\sum_{n=0}^\infty\left( \ell^{2-2\varepsilon} c(n\ell^2+\Delta)
+ \left(\frac{\theta}{\ell}\right)
\ell^{\frac{1}{2}-\varepsilon}\left(\frac{n-\Delta}{\ell}\right)c(n) 
+
c\left(\frac{n-\Delta}{\ell^2}\right)\right)q^n,
\end{equation*}
where
$$
\theta = (-1)^{\frac{1}{2}-\frac{1}{2}({r}+s)}2p^{s}.            
$$
\end{theorem}

\begin{theorem}[{ Newman \cite[Lemma 7 and Theorem 1]{Ne1962} }]
\mylabel{thm:newmanG}
    At $i \infty$, $G(\tau)$ has a pole of order
    $\lfloor\frac{\Delta}{\ell^2}\rfloor$ at most if $\Delta \geq 0$
    and a pole of order $-\Delta$ if $\Delta < 0.$ At 0, $G(\tau)$
    has a pole of order $\lfloor\frac{\Delta^*}{\ell^2}\rfloor$
    at most if $\Delta^* \geq 0$ and a pole of order $-\Delta^*$
    if $\Delta^* < 0.$ Additionally, $G(\tau)$ satisfies
    \begin{equation}
    \mylabel{eq:Gtrans}
       G\left(-\frac{1}{p\tau}\right) = G^*(\tau).
    \end{equation}
\end{theorem}
In this theorem orders of poles at cusps are according to corresponding
uniformizing variables. For $\Gamma_0(p)$ the inequivalent cusps are
$0$ and $i\infty$. The order of $G$ at $0$ corresponds to the
order of $G^{*}$ at $i\infty$.

\section{Proofs} 
\mylabel{sec:proofs}
In this section we prove Theorems \thm{gencong} and \thm{alcong212}.
Suppose $p=2$, $3$, $5$, $7$ or $13$. Then $\Gamma_0(p)$ has
genus $0$ and the function
\beq
\mylabel{eq:Phipdef}
\Phi_p(\tau) = \left( \frac{ \eta(p\tau) }{ \eta(\tau) }\right)^k,
\eeq
where $k = \frac{24}{p-1}$, is a Hauptmodul for $\Gamma_0(p)$.
See Maier \cite{Ma09} and Apostol \cite[pp.86-88]{Ap-BOOK1990}. 
We have the following table for reference:
$$
\begin{array}{c|c c c c c}
p & 2 & 3 & 5 & 7 & 13 \\
\hline
k=\tfrac{24}{p-1} & 24 & 12 & 6 & 4 & 2
\end{array}  
$$
We recall
$$
\eta\PF{-1}{\tau} = \sqrt{-i \tau } \, \eta(\tau),
\qquad \mbox{(see Koblitz \cite[Proposition 14, p.121]{Koblitz-BOOK1993})}
$$
so that 
\beq
\mylabel{eq:Phiptrans}
\Phi_p\PF{-1}{p\tau} = \Parans{\frac{\eta\PF{-1}{\tau}}{\eta\PF{-1}{p\tau}}}^k
= \frac{1}{p^{k/2} \, \Phi_p(\tau) }.
\eeq

Now we suppose that $s>0$ and $r<0$ are integers of opposite parity 
satisfying
$$
0< r + sp < 24 \qquad \mbox{and} \qquad s + rp < 0.
$$
We follow the steps of Section \subsect{newmanres} and recall the functions
\begin{flalign*}
B(\tau) &= \eta^r(\tau) \, \eta^s(p \tau), & 
B^*(\tau) &= \eta^s(\tau) \, \eta^r(p \tau), \\
\phi(\tau) &= 
\prod_{n=1}^\infty (1 - q^n)^r \prod_{n=1}^\infty (1 - q^{pn})^s  &
\phi^*(\tau) &= 
\prod_{n=1}^\infty (1 - q^n)^s \prod_{n=1}^\infty (1 - q^{pn})^r \\
&=\sum_{n=0}^\infty c(n) q^n, & &=\sum_{n=0}^\infty c^*(n) q^n, \\     
g(\tau) &= \frac{B(\ell^2 \tau)}{B(\tau)}, &
g^*(\tau) &= \frac{B^*(\ell^2 \tau)}{B^*(\tau)},
\end{flalign*}
where $3 < \ell \ne p $ is prime.
By Theorem \thm{Gthm},
\begin{align*}
G(\tau) &= 
\Tr{\Gamma_0(\ell^2 p)}{\Gamma_0(p)}{g} \\
&= \frac{1}{\phi(\tau)}
\sum_{n=0}^\infty\left( \ell^{2-2\varepsilon} c(n\ell^2+\Delta)
+ \left(\frac{\theta}{\ell}\right)
\ell^{\frac{1}{2}-\varepsilon}\left(\frac{n-\Delta}{\ell}\right)c(n) 
+
c\left(\frac{n-\Delta}{\ell^2}\right)\right)q^n
\end{align*}
is a weakly holomorphic modular function on $\Gamma_0(p)$.
From Theorem \thm{newmanG}, $G(\tau)$ is holomorphic at $i\infty$
since
$$
0 \le \Delta = (r + sp) \frac{1}{24}(\ell^2-1) < \ell^2 -1
$$
and
$$
\left\lfloor \frac{\Delta}{\ell^2} \right\rfloor = 0.
$$
Also $G(\tau)$ has a pole at $0$ of order $-\Delta^*$ since
$$
\Delta^* = (s + rp) \frac{1}{24}(\ell^2-1) < 0.
$$     
We have 
\begin{align}
\mylabel{eq:Gsid}
G^*(\tau) &= 
\Tr{\Gamma_0(\ell^2 p)}{\Gamma_0(p)}{g^*} \\
&= \frac{1}{\phi^*(\tau)}
\sum_{n=\Delta^*}^\infty\left( \ell^{2-2\varepsilon} c^*(n\ell^2+\Delta^*)
+ \left(\frac{\theta^*}{\ell}\right)
\ell^{\frac{1}{2}-\varepsilon}\left(\frac{n-\Delta^*}{\ell}\right)c^*(n) 
+
c^*\left(\frac{n-\Delta^*}{\ell^2}\right)\right)q^n
\nonumber
\end{align}
is a weakly holomorphic modular function on $\Gamma_0(p)$ 
where 
$$
\theta^* = (-1)^{\frac{1}{2}-\frac{1}{2}({r}+s)}2p^{r}.            
$$
The coefficients in the $q$-expansions of $\phi(\tau)$ and $\phi^*(\tau)$
are integers and clearly the coefficients of $G^*(\tau)$ are
$p$-integral rationals (the only possible denominators are powers
of $\ell$). The function $G^*(\tau)$ has a pole of order
$d^* = - \Delta^*$ at $i\infty$ and is holomorphic at $0$.
It follows that there are $p$-integral rational numbers $\gamma(j)$ such
that
\beq
\mylabel{eq:Gspolys}
G^*(\tau) 
= \lambda_{r,s,p,\ell} + \sum_{j=1}^{d^*} \gamma(j) \Phi_p(\tau)^{-j}
\eeq
for some $p$-integral rational constant $\lambda_{r,s,p,\ell}$.
By Theorem \thm{newmanG} and equation \myeqn{Phiptrans} we have
\beq
\mylabel{eq:Gpolys}
G(\tau) = G^{*}\PF{-1}{p\tau} 
= \lambda_{r,s,p,\ell} + \sum_{j=1}^{d^*} \gamma(j) p^{kj/2} \, 
\Phi_p(\tau)^j.
\eeq
We obtain the result \eqn{gencong} by multiplying by $\phi(\tau)$
and reducing mod $p^{\frac{k}{2}}$. This completes the proof of
Theorem \thm{gencong}.

We can apply the Theorem with $r=-2$, $s=3$, and $p=2$ since
$$
0 < r + sp = 4 < 24,\qquad\mbox{and}\qquad s+rp = -1 < 0.
$$
We find that  
$$
\varepsilon=\frac{1}{2},\quad \theta=16, \quad 
\Delta= \frac{1}{6}(\ell^2-1)=
N_\ell,\quad k=24.
$$
Theorem \thm{alcong212} follows easily. Here we note, in this case,
 that all the
coefficients of $G(\tau)$ and $G^*(\tau)$ are integers.

\section{Faber Polynomials and Further Congruences}
\mylabel{sec:fabpoly}
The Faber polynomials \cite{Fa1903}, \cite{As-Ka-Ni1997} 
have generating function
\begin{align*}
\sum_{m=0}^\infty J_m(x) q^m &= \frac{ E_4(z)^2 E_6(z)}{\Delta(z)} 
\cdot \frac{1}{j(z)-x} \\
&= 1 + (x-744)q + (x^2 -1488 x + 159768) q^2 + \cdots,
\end{align*}
where $j(z)$ is Klein's modular invariant, $E_4(z)$, $E_6()$ are
Eisenstein series of weights $4$, $6$ respectively, and
$$
\Delta(z) = \eta(z)^{24} = q \prod_{n=1}^\infty (1-q^n)^{24}
= \sum_{n=1}^\infty \tau(n) q^n,
$$
is the discriminant modular form and generating function for Ramanujan's
tau-function. The functions $j_m(z) = J_m(j(z))$ are connected
to the denominator formula for the Monster Lie algebra.

This paper grew out of a result of Atkin \cite{At1968} who found
the following theorem for the partition function, $p(n)$, and used
it to prove multiplicative congruences.

\begin{theorem}[Atkin \cite{At1968}]
\mylabel{thm:atkinmod}
Let \begin{equation*}
        Z_{\ell}(\tau) = \sum_{n = -s_{\ell}}^{\infty}
        \left((\ell^3 p(\ell^2n - s_{\ell}) +
        \ell\left(\frac{12}{\ell}\right)\left(\frac{1-24n}{\ell}\right)p(n)
        + p\left(\frac{n+s_{\ell}}{\ell^2}\right)\right)q^{n}
        \mylabel{eq:atkinmod}
    \end{equation*}

where $s_{\ell} = \frac{\ell^2-1}{24}$ and $\ell > 3$ prime. Then
$Z_{\ell}(\tau) (q)_{\infty}$ is a weakly holomorphic modular
function on $\Gamma(1)$, and thus it is a polynomial in $j(\tau)$.
\end{theorem}
 Ono \cite{On2011} connects the polynomials in Theorem \thm{atkinmod}
 with Faber polynomials.
\begin{theorem}[Ono \cite{On2011}]
\mylabel{thm:Onothm}
If $\ell \geq 5$ is prime, and $\delta_{\ell} = \frac{\ell^2-1}{24}$,
then
\begin{equation*}
\widetilde{Z}_{\ell}(\tau)(q)_{\infty}
=  \left(\ell\left(\frac{3}{\ell}\right) +
A(\delta_{\ell};j(\tau))\right)
\end{equation*}
where the polynomials $A(m,x)$ are defined by
\begin{equation*}
    \sum_{m=0}^{\infty} A(m,x)q^m =
    (q)_{\infty}\frac{E_4^2E_6}{\Delta}\frac{1}{j(\tau) - x}
\end{equation*}
\begin{equation*}
    = 1 + (x-745)q + (x^2-1489x+160511)q^2+...
\end{equation*}
\end{theorem}
\begin{example}
\begin{align*}
        Z_5(\tau)(q)_{\infty} &= (j(\tau) - 745) - 5
         =  j - 750 \\
        Z_7(\tau)(q)_{\infty}& = (j^2(\tau) - 1489j(\tau) + 160511)
        - 7= 160504-1489j+j^2\\
        Z_{11}(\tau)(q)_{\infty} &= (j^5(\tau) -
        3721j^4(\tau)+...+(-1971682051554)) + 11\\
        &=-1971682051548+247243785602j-2031082648j^2+4553915j^3-3721j^4+j^5.
\end{align*}
\end{example}
Carney, Etropolski, and  Pitman \cite{Ca-Et-Pi2012}  extend Ono's result to powers
of the eta-function. We extend these results to the eta-quotients
connected with Theorem \thm{gencong}.
We have a generating function for the polynomials in $\Phi_p(\tau)^{-1}$
occurring in \myeqn{Gspolys}.
This is a new result, was not previously discovered
by Newman. For $p = 2, 3,5,7,$ and $13$, and $r$ and $s$ satisfying the
conditions of Theorem \thm{gencong}, define polynomials
$\sB_{r,s,p}(m,x)$ by
\begin{equation}
\mylabel{eq:sBgen}
    \sum_{m=0}^{\infty}\sB_{r,s,p}(m,x)q^m =
    \frac{\delta_q(\Phi_p(\tau)^{-1})}{\phi^*(\tau)}
   \cdot \frac{1}{x-\Phi_p(\tau)^{-1}}
\end{equation}
where $\Phi_p(\tau)$ is defined in \myeqn{Phipdef},
and $\delta_q = q\frac{d}{dq}$.
\begin{theorem} 
\mylabel{thm:onoanalog}
Suppose $p=2$, $3$, $5$, $7$ or $13$ and $3<\ell\ne p$ is prime.
Suppose that $s>0$ and $r<0$ are integers of opposite parity satisfying
$$
0< r + sp < 24 \qquad \mbox{and} \qquad -24 < s + rp < 0. 
$$
Then,
  \begin{equation}
\mylabel{eq:onoanalog}
    G^*(\tau) =  \sB_{r,s,p}(-\Delta^*, \Phi_p(\tau)^{-1})
    + \left(\frac{\theta^*}{\ell}\right)
    \ell^{\frac{1}{2}-\varepsilon}\left(\frac{-\Delta^*}{\ell}\right).
\end{equation}
\end{theorem}
\begin{remark}
We have used the notation $\sB_{r,s,p}$ to avoid
confusion with $B_{r,s,p}(\tau)$ defined in \myeqn{Bdef}.
\end{remark}
\begin{proof}
From \myeqn{Gsid}
\begin{align*}
\phi^*(\tau) \, G^*(\tau) &= 
\sum_{n=\Delta^*}^\infty\left( \ell^{2-2\varepsilon} c^*(n\ell^2+\Delta^*)
+ \left(\frac{\theta^*}{\ell}\right)
\ell^{\frac{1}{2}-\varepsilon}\left(\frac{n-\Delta^*}{\ell}\right)c^*(n) 
+
c^*\left(\frac{n-\Delta^*}{\ell^2}\right)\right)q^n\\
&= q^{\Delta^*} + 
\ell^{\frac{1}{2}-\varepsilon}\left(\frac{-\Delta^*}{\ell}\right) +
O(q),
\end{align*}
since
$$
-24 < s + rp < 0
$$
implies
$$
-\ell^2 + 1 < \Delta^* < 0,\qquad\mbox{and}\qquad 1 \le \Delta + \ell^2
$$
and $q^{\Delta^*}$ is the only negative power of $q$ that appears.
Here
$$
\frac{1}{\phi^*(\tau)} = E(q)^{-s} E(q^p)^{-r} 
= \sum_{n=0}^\infty f(n) q^n,
$$
with$f(0)=1$, so that
\begin{align*}
G^{*}(\tau) &= \left(
             q^{\Delta^*} + 
             \ell^{\frac{1}{2}-\varepsilon}\left(\frac{-\Delta^*}{\ell}\right) +
            O(q) \right) \times
       \left( f(0) + f(1) q + f(2) q^2 + \cdots\right)\\
&=
\sum_{j=0}^{-\Delta^*} f(j) q^{\Delta^* _ j} + 
             \ell^{\frac{1}{2}-\varepsilon}\left(\frac{-\Delta^*}{\ell}\right) +
            O(q).                        
\end{align*}
From Beneish and Larson \cite{Be-La2015}, for each $m\ge1$ 
there is a polynomial $P_m(x)$ of degree $m$ such that
\begin{equation*}
    P_m(\Phi_p(\tau)^{-1}) = q^{-m} + O(q)
\end{equation*}
is the unique weakly holomorphic modular function on $\Gamma_0(p)$
of this form whose only pole is at $i\infty$. The polynomials $P_m(x)$
are given by
\begin{equation*}
    \sum_{m=0}^{\infty}P_m(x)q^m  =
    \frac{\delta_q(\Phi_p(\tau)^{-1})}{x-\Phi_p(\tau)^{-1}}
\end{equation*} where $\delta_q = q\frac{d}{dq}.$ Consider
\begin{equation*}
    \sum_{j=0}^{\Delta^*} f(j) P_{(-\Delta^*-j)}(\Phi_p(\tau)^{-1})
    + \left(\frac{\theta^*}{\ell}\right)
    \ell^{\frac{1}{2}-\varepsilon}\left(\frac{-\Delta^*}{\ell}\right)
    = \sum_{j=0}^{\Delta^*} f(j) q^{\Delta^*+j}
    +  \left(\frac{\theta^*}{\ell}\right)
    \ell^{\frac{1}{2}-\varepsilon}\left(\frac{-\Delta^*}{\ell}\right)
    + O(q).
\end{equation*}
It follows that
\begin{equation*}
   G^*(\tau) 
    = \sum_{j=0}^{\Delta^*} f(j) q^{-\Delta^*+j}
    +  \left(\frac{\theta^*}{\ell}\right)
    \ell^{\frac{1}{2}-\varepsilon}\left(\frac{-\Delta^*}{\ell}\right).
\end{equation*}
since the difference has no poles and is zero. Now
\begin{align*}
&    \sum_{m=0}^{\infty}\sB_{r,s,p}(m,x)q^m =
    \frac{1}{\phi^*(\tau)}\frac{\delta_q(\Phi_p(\tau)^{-1})}{x-\Phi_p(\tau)^{-1}}
\\
&\quad
   =
\left(\sum_{m=0}^{\infty}P_m(x)q^m\right)\left(\sum_{n=0}^{\infty}b(n)q^n\right),
\end{align*}
and 
\begin{equation*}
    \sB_{r,s,p}(m,x) = \sum_{j=0}^mb(j) \, P_{m -j}(x),
\end{equation*}
for $m\ge0$.
Thus
\begin{align*}
    G^*(\tau)&=
    \sum_{j=0}^{-\Delta^*}f(j)P_{\Delta^*+j}(\Phi_p(\tau)^{-1})+\left(\frac{\theta^*}{\ell}\right)
    \ell^{\frac{1}{2}-\varepsilon}\left(\frac{-\Delta^*}{\ell}\right)\\
    &= \sB_{r,s,p}(-\Delta^*,\Phi_p(\tau)^{-1})
    + \left(\frac{\theta^*}{\ell}\right)
    \ell^{\frac{1}{2}-\varepsilon}\left(\frac{-\Delta^*}{\ell}\right)
\end{align*}
which proves the theorem.
\end{proof}

Next we obtain a new identity for the eigenvalue
$\lambda_{r,s,p,\ell} $ in Theorem \thm{gencong} when
$r$, $s$, and $p$ satisfy the conditions of Theorem \thm{onoanalog}.
\begin{theorem}
\mylabel{thm:lamdaid}
Suppose $p=2$, $3$, $5$, $7$ or $13$ and $3<\ell\ne p$ is prime.
Suppose that $s>0$ and $r<0$ are integers of opposite parity satisfying
$$
0< r + sp < 24 \qquad \mbox{and} \qquad -24 < s + rp < 0. 
$$
Let
\begin{equation*}
    \frac{E_p^*(q)}{\phi^*(\tau)} = \sum_{n=0}^{\infty} \mu(n)q^n
\end{equation*}
where
\begin{equation*}
    E_p^*(q) = \frac{pE_2(q^p) - E_2(q)}{p-1}.
\end{equation*}
Then
  \begin{equation*}
    \lambda_{r,s,p,\ell} = \ell^{2-2\varepsilon}
    c(\Delta) + \left(\frac{\theta}{\ell}\right)
    \ell^{\frac{1}{2}-\varepsilon}\left(\frac{-\Delta}{\ell}\right)+
    c\left(\frac{-\Delta}{\ell^2}\right) =
    \mu({\Delta^*})+\left(\frac{\theta^*}{\ell}\right)
    \ell^{\frac{1}{2}-\varepsilon}\left(\frac{-\Delta^*}{\ell}\right).
\end{equation*}
\end{theorem}
\begin{remark}
The function
$$
E_2(q) =  1 - 24 \sum_{n=1}^\infty \sigma(n) q^n
$$
is Ramanujan's quasi-modular form of weight $2$.
The coefficients $c(n)$ and $c^*(n)$ are defined in \myeqn{phidef},
\myeqn{phistardef}.
\end{remark}
\begin{proof}
From \myeqn{Gpolys} we have
$$     
G(\tau) 
= \lambda_{r,s,p,\ell} + \sum_{j=1}^{d^*} \gamma(j) p^{kj/2} \, 
\Phi_p(\tau)^j.
$$        
Thus $\lambda_{r,s,p,\ell}$ is the coefficient of $q^0$ in $G(\tau)$.
By Theorem \thm{onoanalog},
   \begin{equation*}
      G^*(\tau) = \sB_{r,s,p}(-\Delta^*, \Phi_p^{-1})
      + \left(\frac{\theta^*}{\ell}\right)
      \ell^{\frac{1}{2}-\varepsilon}\left(\frac{-\Delta^*}{\ell}\right),
   \end{equation*}
so that by \myeqn{Gtrans} we have
\begin{equation*}
    G(\tau) = \sB_{r,s,p}(-\Delta^*, p^{\frac{k}{2}}\Phi_p)
    + \left(\frac{\theta^*}{\ell}\right)
    \ell^{\frac{1}{2}-\varepsilon}\left(\frac{-\Delta^*}{\ell}\right).
\end{equation*}
Comparing coefficients of $q^0$ leads to
\begin{equation*}
    \lambda_{r,s,p,\ell} = \sB_{r,s,p}(-\Delta^*, 0)
    + \left(\frac{\theta^*}{\ell}\right)
    \ell^{\frac{1}{2}-\varepsilon}\left(\frac{-\Delta^*}{\ell}\right)c^*(0).
\end{equation*} 
From \myeqn{sBgen} we see that $\sB_{r,s,p}(-\Delta^*,0)$ is the coefficient
of $q^{-\Delta^*}$ in
\begin{equation*}
    \frac{1}{\phi^*(\tau)}\biggl\{\frac{\delta_q(\Phi_p(\tau)^{-1})}{-\Phi_p(\tau)^{-1}}\biggr\}.
\end{equation*}
It is not difficult 
to show that
\begin{equation*}
    \delta_q(\Phi_p(\tau)^{-1}) = -\Phi_p(\tau)^{-1}\left(\frac{p
    E_2(q^p) - E_2(q)}{p-1}\right),
\end{equation*}
so that
\begin{equation*}
    \frac{1}{\phi^*(\tau)}\biggl\{\frac{\delta_q(\Phi_p(\tau)^{-1})}{-\Phi_p(\tau)^{-1}}\biggr\}
    = \frac{E_p^*(q)}{\phi^*(\tau)} = \sum_{n \geq 0}\mu(n)q^n,
\end{equation*}
and we have $\sB_{r,s,p}(-\Delta^*, 0) = \mu(\Delta^*)$.
Therefore,
\begin{equation*}
    \lambda_{r,s,p,\ell} = \sB_{r,s,p}(-\Delta^*, 0) 
    +\left(\frac{\theta^*}{\ell}\right)
    \ell^{\frac{1}{2}-\varepsilon}\left(\frac{-\Delta^*}{\ell}\right)
    = \mu({-\Delta^*})+\left(\frac{\theta^*}{\ell}\right)
    \ell^{\frac{1}{2}-\varepsilon}\left(\frac{-\Delta^*}{\ell}\right).
\end{equation*}
Therefore by Theorem \thm{gencong} we have 
  \begin{equation*}
    \lambda_{r,s,p,\ell} = \ell^{2-2\varepsilon}
    c(\Delta) + \left(\frac{\theta}{\ell}\right)
    \ell^{\frac{1}{2}-\varepsilon}\left(\frac{-\Delta}{\ell}\right)+
    c\left(\frac{-\Delta}{\ell^2}\right) =
    \mu({\Delta^*})+\left(\frac{\theta^*}{\ell}\right)
    \ell^{\frac{1}{2}-\varepsilon}\left(\frac{-\Delta^*}{\ell}\right).
\end{equation*}
\end{proof}

Next we make explicit Theorem \thm{Gthm} for the case $r=-2$, $s=3$, $p=2$.
In this case
\begin{align*}
\phi(\tau) &= 
\prod_{n=1}^\infty
\frac{ (1-q^{2n})^3 }{ (1-q^{n})^2}
=  \sum_{n=0}^\infty \alpha(n) q^n =
 A(q),\\
\Delta &= \frac{1}{6}(\ell^2 -1) = N_\ell,\quad \varepsilon = \frac{1}{2}, \quad
\theta = 16,\\
\Delta^* &= -\frac{1}{24}(\ell^2 -1) = -s_\ell,\quad \varepsilon^* = \frac{1}{2}, \quad
\theta^* = \frac{1}{2},\\
\phi^*(\tau) &= 
\prod_{n=1}^\infty
\frac{ (1-q^{n})^3 }{ (1-q^{2n})^2}
=  \sum_{n=0}^\infty \beta(n) q^n = B(q).   
\end{align*}
It turns out that
$$
\beta(n) = M_e(n) - M_o(n),
$$
for $n>1$, where $M_e(n)$ (resp. $M_o(n)$) denote the number of 
partitions of $n$ with even (resp. odd) crank. See \cite{Ch-Ch-Ga2024}
for details.

By Theorem \thm{Gthm} we have
\begin{theorem}
\mylabel{thm:Ztilde} 
Let $\ell > 3$ be prime, $N_{\ell} = \frac{1}{6}(\ell^2-1)$, and 
\begin{equation*}
    \widetilde{Z}_{\ell}(\tau) =  \frac{1}{A(q)}\sum_{n
    \geq 0} \left(\ell \alpha(\ell^2n + N_{\ell}) +
    \left(\frac{-1}{\ell}\right)\left(\frac{N_{\ell}-n}{\ell}\right)\alpha(n)
    +
    \alpha\left(\frac{n-N_{\ell}}{\ell^2}\right)\right)q^{n}.
\end{equation*}
Then $\widetilde{Z}_{\ell}$ is a weakly holomorphic modular function
on $\Gamma_0(2)$, with only a pole at $0$ so that it is 
a polynomial in 
$\Phi_2(\tau) = \left( \frac{ \eta(2\tau) }{ \eta(\tau) }\right)^{24}$,
a Hauptmodul for $\Gamma_0(2)$.                          
\end{theorem}
\begin{example} 
    We have
    \begin{equation*}
        \widetilde{Z}_5(\tau) = 26+4096\Phi_2(\tau),
    \end{equation*}
    \begin{equation*}
        \widetilde{Z}_7(\tau) =
        104+208896\Phi_2(\tau)+16777216\Phi_2(\tau)^2,
    \end{equation*}
    \begin{equation*}
        \widetilde{Z}_{11}(\tau) =
        2100+64609\cdot2^{13}\Phi_2(\tau)+2^{27}\cdot19\cdot373\Phi_2(\tau)^2+2^{36}\cdot5^2\cdot11\cdot17\Phi_2(\tau)^2+3\cdot41\cdot2^{48}\Phi_2(\tau)^4+2^{60}\Phi_2(\tau)^5.
    \end{equation*}
    and
    \begin{equation*}
         \widetilde{Z}_{13}(\tau) =
         9294+3^3\cdot5\cdot11^2\cdot199\cdot2^{13}\Phi_2(\tau)+\dots+2^{84}\Phi_2(\tau)^7.
    \end{equation*}
\end{example}

Similarly we have an analogue for the crank parity function.
\begin{theorem}
\mylabel{thm:crankparity}
   Let \begin{equation*}
       B(q) = \frac{E(q)^3}{E(q^2)^2} = \sum_{k\geq 0}\beta(k)q^k
   \end{equation*} where $\beta(k) = M_e(k) - M_o(k)$, and
   \begin{equation*}
       \widetilde{W}_{\ell}(\tau) =  \frac{1}{B(q)}\sum_{n
       \geq -s_\ell} \left(\ell \beta(\ell^2n - s_{\ell}) +
       \left(\frac{3}{\ell}\right)\left(\frac{24n-1}{\ell}\right)\beta(n)
       +
       \beta\left(\frac{n+s_{\ell}}{\ell^2}\right)\right)q^{n},             
   \end{equation*}
    where $s_{\ell} = \frac{1}{24}(\ell^2-1)$. Then,
    $\widetilde{W}_{\ell}(\tau)$ is a weakly holomorphic
    modular function on $\Gamma_0(2)$, with only a pole at $i\infty$ so 
 that it is a polynomial in
    $\Phi_2(\tau)^{-1}$.
\end{theorem}

\begin{example}
We have
\begin{equation*}
    \widetilde{W}_5(\tau) = q^{-1} + 2 + 276q - 2048q^2 + ... =
    26+\Phi_2(\tau)^{-1},
\end{equation*}
\begin{equation*}
    \widetilde{W}_7(\tau) =
    104+51\Phi_2(\tau)^{-1}+\Phi_2(\tau)^{-2},
\end{equation*}
\begin{equation*}
    \widetilde{W}_{11}(\tau) =
    2100+129218\Phi_2(\tau)^{-1}+...+\Phi_2(\tau)^{-5},
\end{equation*}
\begin{equation*}
    \widetilde{W}_{13}(\tau) =
    9294+6501330\Phi_2(\tau)^{-1}+...+\Phi_2(\tau)^{-7}.
\end{equation*}
\end{example}

Theorems \thm{Ztilde}, \thm{crankparity} were independently
discovered by the second author. Many details of proof are missing
from Newman's paper. Full details of relevant transformation
formulas as given in Chapter 3 of the second author's thesis \cite{Mo-thesis}. 
Finally in this section we apply Theorem \thm{lamdaid} 
to determine $c_\ell \pmod{32}$ given in Theorem \thm{alcong212}
in order to obtain Theorem \thm{alcongs32}.
We apply Theorem \thm{lamdaid} with $r=-2$, $s=3$ and $p=2$.
In this case
$$
\sum_{n=0}^\infty \mu(n) q^n = \frac{E_2^{*}(q)}{B(q)}
$$
where
$$
B(q) = \phi^*(\tau) = \frac{E(q)^3}{E(q^2)}.
$$
The function
$$
E_2^*(q) = 2 E_2(q^2) - E_2(q)
$$
is a holomorphic modular form of weight $2$ on $\Gamma_0(2)$.
We note that the eta-quotients
$$
\frac{\eta(\tau)^8}{\eta(2\tau)^4},\qquad 
\frac{\eta(4\tau)^8}{\eta(2\tau)^4} 
$$
are also holomorphic modular forms of weight $2$ on $\Gamma_0(2)$.
See for instance Ono \cite[Theorem 1.64, p.18]{Onbook04}. 
By a standard argument using the
Valence Formula \cite[Theorem 4.1.4, p.98]{Ra}  we easily find that
$$
E_2^*(q) = 
\frac{\eta(\tau)^8}{\eta(2\tau)^4} + 
32\,\frac{\eta(4\tau)^8}{\eta(2\tau)^4}.
$$
From K\"ohler \cite[Eq. (8.19)]{Ko-BOOK2011}, we have
formula:
\begin{equation*}
        \frac{\eta(\tau)^5}{\eta(2\tau)^2} = \sum_{\substack{n > 0 \\
        n \text{ odd}}} \left(\frac{n}{3}\right) nq^{\frac{n^2}{24}},
\end{equation*} so that
\begin{equation*}
        \frac{E(q)^5}{E(q^2)^2} = \sum_{\substack{n > 0 \\ n \text{
        odd}}} \left(\frac{n}{3}\right) nq^{\frac{n^2-1}{24}}.
\end{equation*} 
It follows that 
$$
\mu(s_\ell) \equiv \ell \, \leg{\ell}{3} \pmod{32}.
$$
By Theorem \thm{lamdaid}  
\begin{align*}
c_\ell &= \lambda_{-2,3,2,\ell} =  
    \mu({-\Delta^*})+\left(\frac{\theta^*}{\ell}\right)
    \ell^{\frac{1}{2}-\varepsilon}\left(\frac{-\Delta^*}{\ell}\right)\\
&= \mu(s_\ell) + \leg{2}{\ell}\,\ell\,\leg{s_\ell}{\ell}\\
&= \mu(s_\ell) + \ell\,\leg{-3}{\ell}\\
&\equiv \leg{\ell}{3}( \ell + 1) \pmod{32},
\end{align*}
since $\leg{-3}{\ell}=\leg{\ell}{3}$ by Quadratic Reciprocity.
Theorem \thm{alcongs32} follows easily.

\section{Some congruences for the overpartition function}
  An \textit{overpartition} is a partition in which the first occurrence of
  a part may be overlined. Let $\pp(n)$ denote the number of overpartitions
  of $n$. For example, there are 8 overpartitions of 3:
  $$\begin{aligned} 3,\ {\overline{3}},\ 2+1,\ {\overline{2}}+1,\ 2+{\overline{1}}, \
   {\overline{2}}+{\overline{1}},\ 1+1+1,\ {\overline{1}}+1+1. \end{aligned}$$
  so that $\pp(3)=8$. We have the generating function
  $$
  \overline{P}(q) = \sum_{n=0}^\infty \pp(n) q^n = \frac{(-q;q)_\infty}{(q;q)_\infty}
  = \frac{(q^2;q^2)_\infty}{(q;q)_\infty^2} = \frac{E(q^2)}{E(q)^2}
  = \frac{\eta(2\tau)}{\eta(\tau)^2}.
  $$
  Many authors have proved congruences for $\pp(n)$ modulo powers
  of $2$ starting with Mahlburg \cite{Ma2004}. Recently 
  Yao \cite{Ya2018} has found explicit congruences modulo $2^6$ and $2^{10}$
  and Xue and Yao \cite{Xu-Ya2021} have explicit congruences modulo 
  $2^{11}$.
  Here we apply our Theorems  \thm{gencong}, \thm{Gthm} and \thm{newmanG}
  to obtain identities and congruences modulo $2^{12}$.

  We apply our theorems with $r=-2$, $s=1$, $p=2$ and $\ell$ is any odd prime.
  We note that the results still hold for $\ell=3$ since in
  this case $\Delta^*$ is an integer.
  We find that
  $$
  \varepsilon=-\tfrac{1}{2}, \quad
  k=24, \quad
  \theta=-4, \quad
  \theta^*=-\frac{1}{2}, \quad
  \Delta=0, \quad
  \Delta^{*} = - \frac{1}{8}(\ell^2 - 1).
  $$
  Here
\begin{align*}
  \phi(\tau) &= \overline{P}(q) = \sum_{n=0}^\infty \pp(n) q^n,\\
  \phi^{*}(\tau) &= \overline{Q}(q) 
   = \sum_{n=0}^\infty  \qq(n) q^n
   = \sum_{n=0}^\infty (-1)^n \podd(n) q^n,\\
  &=\frac{(q;q^2)_\infty}{(q;q)_\infty}= \frac{E(q)}{E(q^2)^2} = 
   q^{\tfrac{1}{8}} \frac{\eta(\tau)}{\eta(2\tau)^2},
\end{align*}
where $\podd(n)$ is the number of partitions of $n$ in which only the odd
parts are distinct. Theorems \thm{Gthm} and \thm{newmanG}
imply that
\begin{align*}
\mathcal{Q}_{\ell}(\tau)= G^*(\tau) = \frac{1}{\overline{Q}(q)} 
\sum_{n=-M_\ell}^\infty \left(
    \ell^3 \qq(\ell^2n - M_\ell)
+ \ell\,\leg{1-8n}{\ell}\qq(n) 
  + \qq\PF{n+M_\ell}{\ell^2}\right)  q^n
\end{align*}
is weakly holomorphic modular function on $\Gamma_0(2)$ whose only  pole
is at $i\infty$, where $M_\ell=\frac{1}{8}(\ell^2-1)$, and
$\qq(n) = (-1)^n \podd(n)$. For example, we find that
\begin{align*}
\mathcal{K}_3(\tau) &= 28 + \Phi_2(\tau)^{-1}, \\
\mathcal{K}_5(\tau) &= 126 + 948 \Phi_2(\tau)^{-1}              
                       + 73 \Phi_2(\tau)^{-2} + \Phi_2(\tau)^{-3}.
\end{align*}
Similarly we have that 
\begin{align*}
\mathcal{L}_{\ell}(\tau)= G(\tau) = \frac{1}{\overline{P}(q)} 
\sum_{n=0}^\infty \left(
    \ell^3 \pp(\ell^2n)
+ \ell\,\leg{n}{\ell}\pp(n) 
  + \pp\PF{n}{\ell^2}\right)  q^n
\end{align*}
is weakly holomorphic modular function on $\Gamma_0(2)$ whose only  pole
is at $0$ of order $\frac{1}{8}(\ell^2-1)$.
For example, we find that
\begin{align*}
\mathcal{L}_3(\tau) &= 28 + 2^{12} \, \Phi_2(\tau),\\
\mathcal{L}_5(\tau) &= 
                       126 + 948 \cdot \, 2^{12} \, \Phi_2(\tau)              
                       + 73 \cdot \,2^{24}\, \Phi_2(\tau) + 
                       + 2 ^{36}\, \Phi_2(\tau).
\end{align*}
Applying Theorem \thm{gencong} we have $\lambda_{-2,1,2,\ell} = (\ell^3+1)$
and
\begin{theorem}
\mylabel{thm:Heckeoverptn}
Suppose $\ell$ is any odd prime. Then
    \begin{equation*}
   \ell^3 \overline{p}(\ell^2n) +
   \left(\frac{-n}{\ell}\right)\,\ell\,\overline{p}(n) +
   \overline{p}\left(\frac{n}{\ell^2}\right)q^{n} \equiv (\ell^3 +
   1)\overline{p}(n) \, \, (\bmod{\,2^{12}}).
\end{equation*}
\end{theorem}

From this, we obtain the following 
\begin{cor}       
\mylabel{cor:Heckeoverptn}
    Suppose $\ell$ is an odd prime, 
    and $\ell \equiv -1 \pmod{2^i}$. Then
\begin{equation*}
  \overline{p}(\ell^3n) \equiv 0  \, \, (\bmod{\,2^{i}}).
\end{equation*}
if $\ell\nmid n$.
\end{cor}     

\section{Conclusion and Some Problems}
\mylabel{sec:end}
This paper first arose by considering analogs of Atkin's multiplicative
congruences for the partition function. We observed that both the
crank parity function and the Andrews's $\EObara$-function 
appear in Ramanujan's identities for third order mock theta functions.
See Watson \cite[p.63]{Wa1936}: 
\begin{align*}
f(q) + 4 \phi(-q) &= \vartheta_4(0,q) \prod_{n=1}^\infty  (1 + q^n)^{-1},\\
\nu(q) + q \omega(q^2) &= \tfrac{1}{2} q^{-\tfrac{1}{4}} \vartheta_2(0,q) 
\prod_{n=1}^\infty (1+ q^{2n}).
\end{align*}
The infinite products on the right sides are the generating functions
for the crank parity and $\EObara$-functions:
\begin{align*}
 \vartheta_4(0,q) \prod_{n=1}^\infty  (1 + q^n)^{-1} &=
\frac{ E(q)^3 }{E(q^2)^2 } = \sum_{n=0}^\infty (M_e(n)-M_o(n))q^n,\\
 \tfrac{1}{2} q^{-\tfrac{1}{4}} \vartheta_2(0,q) 
\prod_{n=1}^\infty (1+ q^{2n})
&= 
\frac{ E(q^4)^3 }{E(q^2)^2 } = \sum_{n=0}^\infty \EObar{n}q^n.
\end{align*}
The corresponding eta-quotients are related by the Fricke transformation
$\tau \mapsto \frac{-1}{4\tau}$, which is connected with the transformation
in Theorem \thm{newmanG} with $p=2$. The functions $f(q)$, $\psi(q)$, 
$\nu(q)$ and $\omega(q)$ are third order mock theta functions.
Watson \cite[p.75]{Wa1936}  found a similar transformation connecting $f(q)$
and $\omega(q)$. This was exploited by Zwegers \cite{Zw2001} in his completion
of these third mock theta functions to harmonic Maass forms. In a subsequent
paper we obtain congruences mod powers of $2$ for the coefficients
of $\omega(q)$.

We conclude by posing some problems.
\begin{enumerate}
\item[(1)]
In this paper we have restricted to half-integer weight
eta-quotients $B_{r,s,p}(\tau)$ where $\Gamma_0(p)$ has genus $0$.
We pose the problem for more general eta-quotients.
\item[(2)]
The congruences in Theorem \thm{gencong} arise because the function
$G(\tau)$ only has a pole at $0$. It would be interesting
to consider when $G(\tau)$ has poles at both $i\infty$ and $0$.
\item[(3)]
We obtained Theorem \thm{alcongs32} from the congruence
$$
\mu(s_\ell) \equiv \ell \, \leg{\ell}{3} \pmod{32},
$$
which in turn came from the following:
\begin{align*}
\sum_{n=0}^\infty \mu(n) q^n &= \frac{ E_2^*(q)}{B(q)} 
= \frac{E(q)^5}{E(q^2)^2} + 32\,q\, \frac{E(q^2) E(q^4)^8}{E(q)^3 E(q^2)^4}\\
&\equiv 
         \sum_{\substack{n > 0 \\ n \text{odd}}} 
\leg{n}{3} \, n\, q^{\frac{n^2-1}{24}} \pmod{32}.
\end{align*}
Find analogous results for more general $\mu(n)$ which will
lead to congruences for the eigenvalues $\lambda_{r,s,p,\ell}$
in Theorem \thm{gencong}.
\end{enumerate}

\end{document}